\documentclass{article}

\usepackage{amssymb}

\newcommand{\remarks}{{\bf Remarks:  }}

\newcommand{\dimv}{\underline{\dim}}

\newcommand{\ses}[3]
{\mbox{$0 \rightarrow #1 \rightarrow #2 \rightarrow #3 \rightarrow 0$}}

\newtheorem{theorem}{Theorem}[section]
\newtheorem{lemma}[theorem]{Lemma}
\newtheorem{definition}[theorem]{Definition}
\newtheorem{proposition}[theorem]{Proposition}
\newtheorem{corollary}[theorem]{Corollary}

\newtheorem{conjecture}[theorem]{Conjecture}
\newtheorem{hope}[theorem]{Hope}
\newtheorem{question}[theorem]{Question}

\begin{document}

\parindent0pt

\title{The use of geometric and quantum group techniques for wild quivers}

\author{Markus Reineke\\ Fachbereich Mathematik\\ Universit\"at Wuppertal\\ Gau\ss str. 20\\ D - 42097 Wuppertal, Germany\\ email: reineke@math.uni-wuppertal.de}

\date{}

\maketitle

\section{Introduction}\label{introduction}

In recent years, the representation theory of quivers has enjoyed an enormous impact of new techniques both from algebraic geometry and from quantum group theory. Some examples of this are the realization of quantum groups in terms of the representation theory of quivers via C.~M.~Ringel's Hall algebra approach \cite{Ri1,Gr1}, G.~Lusztig's geometric interpretation of the Hall algebra approach \cite{Lu1,Lu2,Lu3} and H. Nakajima's work on quiver varieties \cite{Nak1, Nak2}.\\[1ex]
It is thus very desirable to study whether the relations to these areas, and applications of the deep methods which they contain, can be used for the development of substantially new results, methods and perspectives for the representation theory of wild quivers.\\[1ex]
Although much substantial work already has been done, it is only fair to say that the representation theory of wild quivers still offers many mysteries. One of the main points which will be made in the following is that in studying the representation theory of wild quivers, one has to focus on finding the right questions. \\[1ex]
This is of course not surprising by considering the essentially different nature of representation-finite and tame algebras compared to wild algebras: the body of knowledge in the representation theory of finite dimensional algebras, aimed towards understanding representations of finite type and tame algebras, has developed many techniques for understanding essentially discrete phenomena (like parametrizing the indecomposable representations, Auslander-Reiten theory, etc.). But the nature of the representation theory of wild quivers is a continuous one from the very beginning, in that indecomposables appear with arbitrary numbers of parameters (see however \cite{Ker} for the study of wild quivers using Auslander-Reiten theory).\\[1ex]
The first aim of the present paper is to review some of the results and techniques which are available. The second, maybe more important, one, is to speculate about further developments, to show up possible new questions, and to propose future research projects.\\[1ex]
The material to be reviewed in the following can be roughly divided into two types: first, work directly dealing with (the geometry of) quivers, for example of W.~Crawley-Boevey, H.~Derksen, V.~G.~Kac, A.~King, L.~Le Bruyn, A.~Schofield, M.~Van den Bergh, J.~Weyman and the author. Second, some material from fields like classical invariant theory, algebraic group theory, quantum group theory, vector bundle theory and mathematical physics is discussed or at least mentioned, since it either has relations to quiver theory, or is expected or hoped to have so. One of the aims in reviewing this material is to draw the reader's attention to some developments which might turn out to be useful for the representation theory of quivers in the future. The choice of results to be reviewed is of course completely subjective, chosen according to the author's taste and (restricted) knowledge.\\[1ex]
The paper is organized as follows: in Section \ref{notation}, some basic notation on representations of quivers, their moduli spaces, and their Hall algebras is reviewed. Section \ref{results} mainly reviews some (well-known) results, which should serve as examples of what kinds of questions about representations of wild quivers one can actually answer. Section \ref{loop} tries to formulate some questions and perspectives about the archetypical example of wildness, namely the multiple loop quiver. In Section \ref{strat}, one of the at present most successful techniques for the geometric study of quivers, the development of good stratifications of varieties of representations, is discussed in general. The main aim is to show the need for a development of a general theory of such stratifications. In Section \ref{moduli}, some first results in the beginning theory of moduli of quiver representations are reviewed. The main point here is to inspire the reader to consider the vast literature on moduli spaces of vector bundles, which might serve as a guideline towards the future development in the quiver case. Finally, Section \ref{noncomm_geom} reviews some recent results, all of which are inspired by assuming a hypothetical ``noncommutative algebraic geometry", and looking for the consequences of this point of view.\\[1ex]
It should be noticed that the above organization of this paper is sometimes artifical, since all the subjects to be discussed are highly interrelated.\\[1ex]
Throughout the paper, certain ``questions" are posed and ``hopes" are expressed, together with many questions and proposals sketched along the text. They are however of very different nature: some of them ask for a study of particular examples, some for new classes of examples and phenomenology, some for old and deep conjectures, and some for possibly long-range research projects. Also they are sometimes only vaguely formulated.\\[1ex]
None of these questions and hopes are to be taken too literally, and they definitely should not be viewed as explicit formulations of research projects. Instead, they should give the reader an idea of (what the author thinks of) what questions one might pose, and how the field might develop in the future. \\[1ex]
Finally, it should be mentioned that the Proceedings of at least three ICRAs contain overviews over related subjects. An overview over the use of homological and Auslander-Reiten techniques is contained in \cite{Ker}. Geometric approaches to the representation theory of algebras are reviewed in \cite{KrO,Bo3}. For proofs of many of the results to be reviewed here, the reader should also consult L.~Le Bruyn's book project \cite{LBB}.\\[2ex]
This paper grew out of a talk given at a meeting on ``Perspectives in the representation theory of finite dimensional algebras" at the University of Bielefeld. I would like to thank the organizers of this meeting, mainly H.~Krause, for the opportunity to talk there, and for the idea of organizing such a meeting dealing solely with potential future developments. I would like to thank C.~M.~Ringel for encouragement to put this talk, despite its speculative nature, into written form. Furthermore, I would like to thank K.~Bongartz, W.~Crawley-Boevey, L.~Le Bruyn, A.~Schofield and M.~Van den Bergh for inspiring discussions about many of the topics to be covered.\\[1ex]
This paper was completed while the author enjoyed a research stay at the University of Antwerp, with the aid of a grant of the European Science Foundation in the frame of the Priority Programme ``Noncommutative Geometry".

\section{Notation}\label{notation}

\subsection{General notation}\label{general}

We start by fixing some notation which is used constantly throughout this paper.\\[1ex]
Let $Q$ be a finite quiver with (finite) set of vertices $I$ and (finitely many) arrows $\alpha:i\rightarrow j$ for $i,j\in I$. We will assume throughout that $Q$ does not contain oriented cycles.\\[1ex]
The free abelian group ${\bf Z}I$ of $I$ has canonical basis elements $i$ for $i\in I$. Let ${\bf N}I$ be the subsemigroup of nonnegative linear combinations $d=\sum_{i\in I}d_ii$ of the $i\in I$. The Euler form $\langle\_\, ,\_\, \rangle $ on ${\bf Z}I$ is defined by
$$\langle i,j\rangle=\delta_{i,j}-|\{\mbox{arrows from $i$ to $j$ in $Q$}\}|.$$
For a field $k$, the category of finite dimensional $k$-representations of $Q$ is denoted by $\bmod kQ$. For all general facts on the representation theory of quivers, the reader is referred to \cite{ARS}.\\[1ex]
Given a dimension vector $d\in{\bf N}I$, the variety $R_d$ of $k$-representations of $Q$ of dimension vector $d$ is defined as the affine $k$-space
$$R_d=\bigoplus_{\alpha:i\rightarrow j}{\rm Hom}_k(k^{d_i},k^{d_j}),$$
on which the algebraic group
$$G_d=\prod_{i\in I}{\rm GL}_{d_i}(k)$$
acts via
$$(g_i)_i\cdot(M_\alpha)_\alpha=(g_jM_\alpha g_i^{-1})_{\alpha:i\rightarrow j}.$$
By definition, the orbits ${\cal O}_M$ are in bijection with the isomorphism classes $[M]$ of $k$-representations of $Q$ of dimension vector $d$.

\subsection{Moduli spaces}\label{notation_moduli}

We next define (semi-) stability and moduli spaces for quiver representations. For this, we first have to make a choice. The space of ${\bf Z}$-linear functions ${\rm Hom}_{\bf Z}({\bf Z}I,{\bf Z})$ has a basis consisting of the $i^*$ for $i\in I$, defined by $i^*(j)=\delta_{i,j}$ for $j\in I$. We denote by $\dim$ the function $\sum_{i\in I}i^*$. Let us choose a function
$\Theta\in{\rm Hom}_{\bf Z}({\bf Z}I,{\bf Z})$ , and define the slope function $\mu:{\bf N}I\setminus 0\rightarrow {\bf Q}$ via
$$\mu(d):=\frac{\Theta(d)}{\dim d}.$$
For a representation $X$, we will shortly write $\mu(X)$ for $\mu(\dimv X)$.\\[1ex]
Call a representation $X$ semistable (resp.~stable) if for all proper non-zero subrepresentations $U$ of $X$, we have $\mu(U)\leq\mu(X)$ (resp.~$\mu(U)<\mu(X))$.\\
Note that this notion of semistability is essentially equivalent to that of A.~King from \cite{Ki1}, where a representation $X$ such that $\Theta(\dimv X)=0$ is defined to be semistable if and only if $\Theta(\dimv U)\geq 0$ for all proper non-zero subrepresentations $U$ of $X$. This can be seen as follows: fix the dimension vector $d$, and define a functional
$$\widetilde{\Theta}=\mu(d)\cdot\dim-\Theta.$$
Then it is easy to see from the definitions that a representation $X$ of dimension vector $d$ is semistable in the sense defined above if and only if it is $\widetilde{\Theta}$-semistable in the sense of A.~King.\\[1ex]
The present definition, coming from vector bundle theory (to be discussed in Section \ref{vb}) has the advantage of allowing to define (semi-)stability for representations of arbitrary dimension vectors.\\[1ex]
Denote by $R_d^{ss}$ (resp.~$R_d^s$) the set of points of $R_d$ corresponding to semistable (resp.~stable) representations. By applying the general machinery of D.~Mumford's Geometric Invariant Theory \cite{Mu}, A.~King proved the following \cite{Ki1}:
\begin{theorem}[A.~King] The variety $R_d^{ss}$ is an open subset of $R_d$, and $R_d^s$ is an open subset of $R_d^{ss}$. There exists an algebraic quotient ${\cal M}_d^{ss}(Q):=R_d^{ss}//G_d$
of $R_d^{ss}$ by $G_d$. Moreover, there exists a geometric quotient
${\cal M}_d^s(Q):=R_d^s/G_d$
of $R_d^s$ by $G_d$. The variety ${\cal M}_d^{ss}$ is projective, and ${\cal M}_d^s$ is a smooth open subvariety of ${\cal M}_d^{ss}$.
\end{theorem}

\remarks
\begin{itemize}

\item For precise definitions of the notions of algebraic and geometric quotients, a detailed discussion of these various quotient constructions, with particular emphasis on their applications to representation theory, see \cite{Bo3}.

\item The moduli space ${\cal M}_d^{ss}$ is most easily defined as follows: define a character $\chi$ on $G_d$ via
$$\chi((g_i)_{i\in I})=\prod_{i\in I}\det(g_i)^{\widetilde{\Theta}(i)},$$
where $\widetilde{\Theta}$ refers to the above modification of $\Theta$.
With respect to this character, we can speak about semi-invariants of the action of $G_d$ on $R_d$ of weight $\chi^n$, namely, $k[R_d]^{G_d,\chi^n}$ is defined as the space of regular functions $f:R_d\rightarrow k$ such that $f(gx)=\chi(g)^n\cdot f(x)$ for all $g\in G_d$ and all $x\in R_d$. Then $${\cal M}_d^{ss}={\rm\bf Proj}\bigoplus_{n\geq 0}k[R_d]^{G_d,\chi^n},$$
the projective spectrum of the graded ring of semi-invariants of weight a multiple of $\chi$. This relates the study of quiver moduli to the study of semi-invariants of \cite{DW,DZ,SVDB}.

\item If stability and semistability coincide for representations of dimension vector $d$, then ${\cal M}_d^s={\cal M}_d^{ss}$ is a smooth projective variety, which makes it accessible to many
techniques of algebraic geometry. A priori, this is fulfilled if $\mu(e)\not=\mu(d)$ for all $0\not=e\leq d$. A more handy criterion is the following \cite{Ki1}:
\begin{center}
For generic $\Theta$, we have $R_d^s=R_d^{ss}$ if and only if $d$ is indivisible.
\end{center}

\item The points of ${\cal M}_d^{ss}$ do not in general parametrize the isomorphism classes of semistable representations, but only certain equivalence classes of them. Namely, for each $\mu\in{\bf Q}$, we have the full subcategory $\bmod_\mu kQ$ consisting of semistable representations of slope $\mu$, in which the simple objects are precisely the stable representations of slope $\mu$. Two representations $X,Y\in\bmod_\mu kQ$ are called $S$-equivalent if they have the same Jordan-H\"older type in $\bmod_\mu kQ$. More explicitely, $X$ and $Y$ are $S$-equivalent if they can be filtered by the same finite set of stable representations of slope $\mu$.

\end{itemize}

\subsection{Hall algebras}\label{notation_hall}

In this section, we assume that $k$ is a finite field, and we let $v\in{\bf C}$ be a square root of
its 
cardinality (this unusual choice will become clear in Theorem
\ref{ringel} 
below).
To define the Hall algebra, we use M.~Kapranov's version \cite{Kap1} of 
C.~M.~Ringel's definition \cite{Ri1}, since
it makes the relation to the geometry of the representation spaces
clear.
We define a convolution on functions on representation spaces. Let ${\bf 
C}^{G_d}(R_d)$ be the (finite-dimensional)
complex vector space of $G_d$-invariant functions on $R_d$.
The elements of ${\bf C}^{G_d}(R_d)$ can thus also be viewed as formal
${\bf 
C}$-linear
combinations of the orbits, and therefore, of the isomorphism classes of 
$d$-dimensional representations by the result above. Let
$$H_v(Q)=\bigoplus_{d\in{\bf N}I}{\bf C}^{G_d}(R_d)$$
be the direct sum of all these spaces, which is viewed as an ${\bf
N}I$-graded 
complex vector space in the obvious way.
To define a graded multiplication
on $H_v(Q)$, let $f\in{\bf C}^{G_d}(R_d)$ and $g\in {\bf C}^{G_e}(R_e)$ be
two 
functions. We define the convolution product
$(f*g)\in{\bf C}^{G_{d+e}}(R_{d+e})$ by
$$(f*g)(X)=v^{\langle d,e\rangle}\sum_{U\subset X}f(X/U)g(U),$$
where the sum runs over the (finitely many) subrepresentations $U$ of $X$ of 
dimension
vector 
$e$. Note that the value of $f$ (resp.~$g$)
at $X/U$ (resp.~at $U$) is well-defined, since the values of these
functions depend only on the isomorphism classes of
representations.
It is then easy to see that $H_v(Q)$ becomes an associative ${\bf
N}I$-graded 
${\bf C}$-algebra.\\[1ex]
For a vertex $i\in I$, the variety $R_i$ consists of single point;
let $E_i$ be the function having value $1$ on this point. The subalgebra 
$C_v(Q)$ of $H_v(Q)$ generated by the elements $E_i$ for
$i\in I$ is
called the composition algebra. To formulate the main result about these 
constructions, we recall some notation on quantized enveloping
algebras.\\[2ex]
In the theory of Kac-Moody Lie algebras \cite{Kac2}, a (symmetric)
generalized 
Cartan matrix is associated to any finite graph.
Applying this
to the unoriented graph underlying $Q$, we arrive at the Cartan matrix 
$C=(a_{ij})_{i,j\in I}$ defined by $a_{ii}=2$ for all
$i\in I$,
and $-a_{ij}$
equals the number of arrows between $i$ and $j$ in $Q$ (in either
direction) 
for $i\not=j$. Note that $C$ is the symmetrization of the
matrix representing
the bilinear form $\langle\_\, ,\_\, \rangle$. 
Associated to $C$, we have (\cite{Ja1} and \cite{Lu3}) the
positive part 
${\cal U}^+_v(\mathfrak{g})$ of the quantized
enveloping algebra ${\cal U}_v(\mathfrak{g})$, which is defined as the
${\bf 
Q}(v)$-algebra with generators $E_i$ for $i\in I$ and
the following defining
relations (the quantized Serre relations):
$$\sum_{p+p'=1-a_{ij}}(-1)^{p'}E_i^{(p)}E_jE_i^{(p')}=0\mbox{ for 
}i\not=j\mbox{ in }I.$$
In this definition, the element $E_i^{(n)}$ denotes the divided power 
$([n]!)^{-1´}E_i^n$ for $i\in I$ and $n\in{\bf N}$, where 
$[n]!=[1]\cdot[2]\cdot\ldots\cdot[n]$ is the quantum factorial
defined 
via the quantum number $[n]=(v^n-v^{-n})/(v-v^{-1})\in{\bf Z}[v,v^{-1}]$.
Consider also the ${\bf Z}[v,v^{-1}]$-subalgebra $U^+$ of ${\cal 
U}^+_v(\mathfrak{g})$ generated by all $E_i^{(n)}$ for $i\in I$,
$n\in{\bf 
N}$. This algebra can be specialized to any non-zero complex number $v$
in the 
obvious way, yielding a ${\bf C}$-algebra $U^+_v={\bf C}\otimes_{{\bf 
Z}[v,v^{-1}]}U^+$. In particular, the specialization $U^+_1$ is
isomorphic to 
the enveloping algebra ${\cal U}(\mathfrak{g}^+)$ of the positive part
of the 
Kac-Moody Lie algebra $\mathfrak{g}$.
 
The relation between the Hall algebra and the quantized enveloping
algebra is 
given in the following theorem (see \cite{Gr1} and \cite{Ri1}):
 
\begin{theorem}\label{ringel} The map $\eta:E_i\rightarrow E_i$ extends
to an 
isomorphism of ${\bf N}I$-graded ${\bf C}$-algebras
$$\eta:U^+_v\stackrel{\sim}{\rightarrow}C_v(Q).$$
\end{theorem}

\section{Some known results}\label{results}

\subsection{General representations}\label{general_rep}

As a first step in the geometric study of the representation varieties $R_d$, one can ask for their generic properties; that is, properties of the representations $M$ in some open subset of $R_d$. The idea behind this approach is that much of the complexity in studying the orbit structure of a variety with a group action comes from boundary phenomena, which are excluded if one restricts to sufficiently small open subsets. As a prototype of this phenomenon, consider conjugacy classes of matrices; that is, representation varieties for the quiver consisting of one loop. Then it is well known that the diagonalizable matrices form a dense subset. In a representation theoretic formulation, this means that a general representation is a direct sum of one-dimensional representations of the path algebra $k[t]$. On the other hand, we have the full theory of the Jordan canonical form as a boundary phenomenon, which we avoided by restricting to general representations. The study of general representations for arbitrary quivers was started by V.~G.~Kac \cite{Kac1,Kac15} and continued by A.~Schofield \cite{Sc1}. We review some of their results.\\[2ex]
Given dimension vectors $d,e\in{\bf N}I$, we define $hom(d,e)$ (resp.~$ext(d,e)$) as the minimum value of $\dim{\rm Hom}_Q(M,N)$ (resp.~of $\dim{\rm Ext}^1_Q(M,N)$) for representations $M\in R_d$, $N\in R_e$. Since the functions $\dim{\rm Hom}_Q(\_,\_)$ and $\dim{\rm Ext}^1_Q(\_,\_)$ are upper semicontinuous on $R_d\times R_e$, the minimum is obtained on some open subset of $R_d\times R_e$, and thus $hom(d,e)$ (resp.~$ext(d,e)$) can be viewed as the generic value of the corresponding dimensions.\\[1ex]
Similarly, we define the generic decomposition of the dimension vector $d$ as the decomposition $d=d^1+\ldots+d^n$ such that on an open subset $U$ of $R_d$, all representations $M$ belonging to $U$ decompose into a direct sum $M=M_1\oplus\ldots\oplus M_n$ of indecomposables of dimension vectors $\dimv M_i=d^i$ for $i=1\ldots n$. Obviously, this concept is the most interesting from the point of view of the representation theory of finite dimensional algebras.

\begin{theorem}[Schofield] $ext(d,e)$ equals the maximum of the values $-\langle d',e\rangle$, taken over all dimension vectors $d'\leq d$ such that $ext(d',e-d')=0$.
\end{theorem}

\begin{theorem}[Schofield] $d$ is a Schur root if and only if $\langle e,d\rangle\leq \langle d,e\rangle$ for all dimension vectors $e\leq d$ such that $ext(e,d-e)=0$.
\end{theorem}

\begin{theorem}[Kac] $d=d^1+\ldots+d^n$ is the generic decomposition of $d$ if and only if the following holds:
\begin{enumerate}
\item all $d^i$ are Schur roots,
\item $ext(d^i,d^j)$ equals $0$ for all $i\not=j$.
\end{enumerate}
\end{theorem}

From this theorem we see that one gets a recursive method for determining the generic decomposition, which is a priori not very effective algorithmically. Better algorithms can be found in \cite{Sco2} and in \cite{DW2}.\\[1ex]
Of course it would be very desirable to have an explicit formula instead of an algorithm only, both for practical and for theoretical considerations. A similar phenomenon appears in the description of the cohomology of quiver moduli, see Section \ref{rec}. In this case, the recursion can be resolved by combinatorial means. Moreover, one can expect such recursive phenomena whenever certain types of stratifications of representation varieties are involved, such that the geometric structure of arbitrary strata can be reduced to generic strata for smaller dimension types; see \ref{general_strat} for a more precise formulation of this.\\[1ex]
Motivated by the success for the Harder-Narasimhan recursion and by the improved algorithms \cite{Sco2, DW2}, one might propose:

\begin{hope} Several of the recursions occuring in the geometric study of wild quivers can be resolved to explicit formulas.
\end{hope}

\subsection{Local and birational geometry}\label{local_birational}

As the next step in the geometric study of representation varieties, we can ask for local properties of moduli spaces. To do this, it is not neccessary to restrict to the moduli spaces in the sense of A.~King defined above: in general, there are many open subsets $U\subset R_d$ admitting a geometric quotient $U/G_d$. The stable locus $R_d^s$ is just one natural choice (but note that it has certain maximality properties by D.~Mumford's general theory, \cite{Mu}). 
Since we are only interested in local properties, we can always assume $U$ to be sufficiently small.\\[2ex]
The first local property which received a lot of attention is that of rationality, asking whether a
quotient $U/G_d$ looks like an open subvariety of some affine space. This is a long-standing conjecture in the case of the multiple loop quiver (for an overview, see \cite{LB0} or \cite[Chapter 5]{For}), so the best result one can hope for at the moment is a reduction from arbitrary quivers to this case:

\begin{theorem}[\cite{Sco2}] Let $d\in{\bf N}I$ be a Schur root, and let $U\subset R_d$ be an open subset such that a geometric quotient $U/G_d$ exists. Then $U/G_d$ is birationally equivalent to the algebraic quotient $(k^{n\times n})^p//{\rm GL}_n$, where $n$ equals the greatest common divisor of the $d_i$ for $i\in I$, and $p$ equals $1-\langle \frac{d}{n},\frac{d}{n}\rangle$.
\end{theorem}

More generally, we can ask for the geometric structure of an open neighborhood of a given point of a quotient. Again one can get a general result \cite{AL}, but one has to sacrifice the use of the Zariski topology and pass to the etale topology. Alternatively, the reader might think about the base field ${\bf C}$, and asking for open neighborhoods in the analytic topology only. To state the result, we need some notation.\\[1ex]
The description of the points of ${\cal M}_d^{ss}(Q)$ as $S$-equivalence classes of semistables can be rephrased as saying that ${\cal M}_d^{ss}$ parametrizes isomorphism classes of representations which are direct sums of stable representations of the same slope $\mu(d)$. Let $X=X_1^{m_1}\oplus\ldots\oplus X_s^{m_s}$ be such a representation, such that the $X_i$ are pairwise non-isomorphic stables of slope $\mu(d)$. Introduce the so-called local quiver $Q_X$ with set of vertices $\{1,\ldots,s\}$ and $\delta_{i,j}-\langle \dimv X_i,\dimv X_j\rangle$ arrows from $i$ to $j$, for all $i,j=1\ldots s$. Define a dimension vector $d_X$ for $Q_X$ by $d_i=m_i$ for $i=1\ldots s$.

\begin{theorem}[\cite{AL}] In the etale topology, a neighborhood of $X$ in ${\cal M}_d^{ss}(Q)$ is isomorphic to a neighborhood of the zero representation in the affine quotient $R_{d_X}(Q_X)//G_{d_X}$.
\end{theorem}

Note that $Q_X$ usually contains many oriented cycles, so that the affine quotient $R_{d_X}(Q_X)//G_{d_X}$ is highly nontrivial.\\[2ex]
To see some analogy to the result of A.~Schofield mentioned before, look at the case where $d$ is some $n$-th multiple of a Schur root, and where $X$ is the $n$-th power of a stable representation in ${\cal M}_{d/n}^{ss}(Q)$. Then the above theorem tells us that in a neighborhood of $X$, the moduli space looks like the quotient variety of a $(1-\langle \frac{d}{n},\frac{d}{n}\rangle)$-tuple of $n$ by $n$ matrices in a neighborhood of the zero tuple of matrices.\\[2ex]
The two theorems cited in this section indicate a general principle: the study of the structure of the moduli spaces ${\cal M}_d^{ss}$ should be subdivisible into two ``orthogonal" directions: 
\begin{itemize} 
\item the study of ${\cal M}_d^{ss}$ in the case of divisible $d$,
\item the study of the quotient varieties $(k^{n\times n})^m//{\rm GL}_n$.
\end{itemize}

This suggests the following:

\begin{question} How can this principle be put into precise structural statements? For example, can one prove a result which directly relates the variety ${\cal M}_{nd}^{ss}$ to the varieties ${\cal M}_d^{ss}$ and some $(k^{n\times n})^m//{\rm GL}_n$ (or related varieties) in case of indivisible $d$?
\end{question}

Note that this is possible at least birationally by \cite{Sco2}.

\subsection{The Kac conjectures}\label{kac_conj}

Still one of the most fundamental results in the representation theory of wild quivers is the theorem of V.~G.~Kac:

\begin{theorem}[\cite{Kac1}] The dimension vectors $d\in{\bf N}I$ admitting an indecomposable representations of dimension vector $d$ correspond bijectively to the positive roots of the Kac-Moody algebra $\mathfrak{g}_Q$ associated to $Q$.
\end{theorem}

It forms the basis of a number of far-reaching conjectures of V.~G.~Kac \cite{Kac15}, two of which were settled at least in a lot of cases recently \cite{CBV}. 

By Kac's theorem, it is clear to expect a close relation between the representation theory of $Q$ and the structure of the corresponding Kac-Moody algebra ${\mathfrak{g}}_Q$. This relation experienced a far-reaching conceptual explanation by C.~M.~Ringel's Hall algebra construction \cite{Ri1}. To quote M.~van den Bergh, it is surprising that, although Hall algebra theory explained so much about the relation between $\bmod kQ$ and ${\mathfrak{g}}_Q$, the Kac conjectures are still open.\\[2ex]
Consider representations of $Q$ over a finite field with $q$ elements ${\bf F}_q$. A representation is called absolutely indecomposable if it remains indecomposable even after extending scalars to an algebraic closure of ${\bf F}_q$ (see \ref{counting} for similar concepts).

Let $a_d(q)$ denote the number of isomorphism classes of absolutely indecomposable representations of ${\bf F}_qQ$ of dimension vector $d$.

\begin{conjecture}[Kac] The following holds for arbitrary roots $d\in{\bf N}I$:
\begin{enumerate}
\item There exists a polynomial $P_d(t)$ such that $P_d(q)$ equals $a_d(q)$ for all prime powers $q$. This polynomial has nonnegative integer coefficients.
\item The value of $P_d(t)$ at $t=0$ equals the multiplicity of the root $d$ in the Kac-Moody algebra correspoding to $Q$.
\end{enumerate}
\end{conjecture}

\begin{theorem}[\cite{CBV}] The above conjectures are true for indivisible $d$.
\end{theorem}

One of the steps of the proof is a translation of the problem to deformed preprojective varieties. This allows us to apply the results of W.~Crawley-Boevey on these varieties \cite{CB1,CB2,CB3} , and their Hyper-K\"ahler structure introduced by H.~Nakajima \cite{Nak1}. The proofs also involve a study of the arithmetics of the involved varieties, related to the Weil conjectures \cite{De}. For the second conjecture, Harder-Narasimhan techniques are used (see \cite{hnsystem}, \ref{examples_strat}, \ref{cohomology}) to relate to quantum group theoretical results.\\[2ex]
Of course the only hope here could be:

\begin{hope} The Kac conjectures hold in general.
\end{hope}

At the moment, it is unclear which techniques could lead to a proof of the Kac conjectures in full generality. So what this hope mainly expresses is a hope for the development of new techniques.

\section{The multiple loop quiver}\label{loop}

Especially from Section \ref{local_birational} it becomes clear that also in a geometric study of representations of wild quivers, the multiple loop quiver enters naturally. Somehow of course, this is not a surprise, since after all, wildness is just defined by this quiver. It is thus neccessary that new techniques have to be invented for the study of this particular case.\\[1ex]
So denote by $L_m$ the quiver with one vertex and $m$ loops. The representation variety for dimension $n\in{\bf N}$ is the space $(k^{n\times n})^m$ of $m$-tuples of $n\times n$-matrices over $k$, on which the group ${\rm GL}_n$ acts by simultaneous conjugation.\\[2ex]
First of all, one has to ask for new classes of examples. This could either mean to restrict to small $n$, or to restrict the study of the action of ${\rm GL}_n$ to some ``nice" subsets $U\subset(k^{n\times n})^m$. The latter case could mean, for example, to look at simple representations, or at representations with given type of endomorphism ring. It might also turn out to be fruitful to study certain combinatorially defined classes of representations of $L_m$ (like some analogues of tree modules) from a geometric perspective.\\[1ex]
In the past, the studies of the representation varieties for small $n$ mostly dealt with explicit classifications of the orbits, so they do not neccessarily reflect the geometry involved. More precisely, one should ask:

\begin{question} What is the geometric structure of the quotient variety $$V_n^{(m)}=(k^{n\times n})^m//{\rm GL}_n$$ for small $m,n\in{\bf N}$?
\end{question}

Note that this quotient (similar to the case of moduli spaces considered before) parametrizes semisimple representations of dimension $n$ over the free algebra $k\langle x_1,\ldots,x_m\rangle$.\\[1ex]
The following results were pointed out to me by L.~Le Bruyn: Assuming $n\geq 2$, we have that: $V_n^{(m)}$ is smooth only if $(m,n)=(2,2)$; it is a complete intersection only if $(m,n)=(2,2)$ or $=(3,2)$ or $=(2,3)$.\\[1ex]
In Section \ref{exkron22}, we will work out a particular example of moduli spaces for the generalized Kronecker quiver. This example also gives a description for the quotient variety of $m$-tuples of $2$ by $2$ matrices with the aid of a simple reduction technique:\\[1ex]
Denote by $K_m$ the $m$-arrow Kronecker quiver; that is, the quiver with set of vertices $\{i,j\}$, and $m$ arrows from $i$ to $j$. Given a nonnegative integer $n$, there is a map
$$R_n(L_m)\rightarrow R_{ni+nj}(K_{m+1})$$ given by
$$(A_1,\ldots,A_m)\mapsto(id_n,A_1,\ldots,A_m).$$
It is easy to see that this map induces an open embedding (we choose the stability $\Theta=i^*$ for the quiver $K_n$)
$$R_n(L_m)//{\rm GL}_n\rightarrow {\cal M}_{ni+nj}^{ss}(K_{m+1}),$$
such that the isomorphism classes of simple representations are precisely those which map to the classes of stable representations.\\[1ex]
In connection with this simple reduction, there is the following problem. One can stratify $R_{ni+nj}(K_{m+1})$ according to the rank of the matrix corresponding to the first arrow, giving strata $S_r$ for $r=0\ldots n$. Then, as seen above, $S_n$ reduces to $m$-tuples of $n$ by $n$ matrices. As another obvious special case, $S_0$ corresponds to $R_{ni+nj}(K_m)$.

\begin{question} How can one describe the intermediate strata $S_r$ for $r=1\ldots n-1$?
\end{question}

L.~Le Bruyn proposes to extend this technique, so that arbitrary oriented cycles in a quiver can be resolved.\\[3ex]
The second possible new direction for the study of $L_m$ is motivated by quantum group theory; namely, by a result of B.~Sevenhant and M.~van den Bergh \cite{BSVDB}. They prove:

\begin{theorem} Whenever $Q$ is a quiver without oriented cycles, the Hall algebra is isomorphic to the positive part of the quantized enveloping algebra of a Borcherds (generalized Kac-Moody) algebra ${\mathfrak{g}}$.
\end{theorem}

Borcherds algebras differ from Kac-Moody algebras in that imaginary simple roots (even with non-trivial multiplicities) are allowed. For an overview, see \cite{Ray}. An obvious question is:

\begin{question} Is the Hall algebra of $L_m$ isomorphic to the quantized enveloping algebra of the positive part $\mathfrak{g}^+$ of some Borcherds algebra ${\mathfrak{g}}_m$?
\end{question}

This is true at least for $m=0,1$. For $m=0$, the algebra ${\mathfrak{g}}$ is just ${\mathfrak{sl}}_2$, and for $m=1$, the algebra ${\mathfrak{g}}_1^+$ is an abelian Lie algebra in countably many generators, corresponding to the classes of zero matrices of various sizes (see \cite{McD}).\\[1ex]
It is interesting to speculate on the structure of this conjectural algebra ${\mathfrak{g}}_m^+$ for arbitrary $m$. It should be a Borcherds algebra of rank one, say with generator $\alpha$ of the root lattice. Essentially, it should be a free Lie algebra. One can also approach this in the spirit of the Kac conjectures. There should be (imaginary) simple roots for the root spaces $n\alpha$ for all $n\in{\bf N}$, with certain multiplicities $c_n$, which should differ considerably from the root multiplicities $d_n$ themselves. 

\begin{question} Is there a representation theoretic interpretation of these multiplicities $c_n$ and $d_n$? And what is the representation theoretic meaning of the Borcherds denominator formula in this context?
\end{question}

The study of the Hall algebras $H(L_m)$ was started by G.~Lusztig in \cite{Lu}. He proposes to look at the subalgebra generated by the characteristic functions $e_n$ of the zero tuple of matrices of various sizes $n\in{\bf N}$. He proves essentially that this is a free associative algebra in the $e_n$, and proposes to view the monomials in the $e_n$ as a canonical basis of this subalgebra, in the spirit of his general theory of canonical bases of quantum groups. The justification for this, which is the main result of G.~Lusztig on this subalgebra, can be formulated as follows:\\[1ex]
Given a sequence $n_1,\ldots,n_s$ of nonnegative integers summing up to $n$, consider the closed subvariety ${\cal E}$ of $R_n(L_m)$ of tuples of matrices which can be simultaneously conjugated to the upper triangular block form
$$\left[\begin{array}{cccc}0_{n_1}&*&*&*\\ 0&0_{n_2}&*&*\\ &&\ldots&\\ 0&0&0&0_{n_s}\end{array}\right].$$

This can be viewed as the image of the first projection map $\pi$ from the variety $X$ of pairs consisting of a point $(A_1,\ldots,A_m)$ of $R_n(L_m)$ and a flag $k^n=F^0\supset F^1\supset\ldots\supset F^s=0$ of linear subspaces of $k^n$ such that $\dim F^{k-1}/F^k=n_k$ for $k=1\ldots s$, such that $A$ and the flag are compatible, in the sense that $A_k(F^l)\subset F^{l+1}$ for all $k=1\ldots m$ and all $l=1\ldots s$.

\begin{theorem} For $m\geq 2$, the map $\pi:X\rightarrow{\cal E}$ is a small resolution in the sense of intersection cohomology theory \cite[Chapter 8]{CG}.
\end{theorem}

The smallness implies that the pushdown of the constant sheaf on $X$ (in the derived category of constructible sheaves) is perverse (up to some shift). Translated to the language of canonical bases \cite{Lu1}, this means that the monomial $e_{n_1}\ldots e_{n_s}$ in $H(L_m)$ should be viewed as an element of a canonical basis of the Hall algebra.\\[1ex]
However, there is no obvious analogue of a PBW basis! A detailed explanation would lead too far away from the present topic; roughly speaking, a natural orthogonalization procedure which could produce such a PBW basis fails. Therefore, it seems neccessary to produce larger natural subalgebras of $H(L_m)$.

\begin{hope} One can find larger subalgebras of $H(L_m)$ which still allow a reasonable structure theory, in particular with respect to natural bases for them. In particular, there should be small resolutions of certain strata in $(k^{n\times n})^m$ related to analogues of canonical bases, and there should be an analogue of a PBW basis. 
\end{hope}

The third aspect to be mentioned about $L_m$ starts at the very foundations of the representation theory of finite dimensional algebras. Roughly speaking, the wildness of a quiver $Q$ is defined by admitting an embedding of the category of representations of some $L_m$ for $m\geq 2$ into 
$\bmod kQ$. One can ask for geometric analogues of this property. This was developed by K.~Bongartz in \cite{Bo3}. One possible question in this direction is:

\begin{question} What is the nature of maps between representation varieties coming from representation embeddings of the corresponding categories? In particular, are there restrictions on the possible images of such embeddings?
\end{question}

\section{Stratifications}\label{strat}

The technical heart of many proofs of results on the geometry of quiver representations is formed by the construction of stratifications with certain desirable properties. They constitute a general tool for the development of reduction techniques, recursions, and other means of simplifying the study of representation varieties, moduli spaces and other quotients. Because of this importance, they deserve a general study for themselves.

\subsection{General discussion of stratifications}\label{general_strat}

In most generality, the term stratification means that we are looking for finite decompositions
$$R_d=\bigcup_{i\in {\cal I}}S_i$$
of a representation variety $R_d$ into locally closed  subvarieties $S_i$, parametrized by certain combinatorial objects. Then one might ask for example the following questions about this stratification:

\begin{itemize}
\item Is there an efficient way to determine whether $S_i$ is non-empty in terms of the combinatorial object $i$?
\item What are the general geometric properties of the strata $S_i$? For example, are they irreducible? What is their dimension? Are they smooth? Are they rational varieties?
\item Is the closure $\overline{S_i}$ of a stratum again a union of some strata? If this is the case, can one describe the relation $S_j\subset\overline{S_i}$ in terms of the objects $i$ and $j$? 
\item Or can one at least define a partial ordering $i\leq j$ on ${\cal I}$ such that
$$\overline{S_i}\subset\bigcup_{i\leq j}S_j?$$
\item Does there exist an algebraic quotient $S_i//G_d$? Or even a geometric quotient $S_i/G_d$? If this is the case, what can one say about the structure of these quotients?
\item Answer the previous question at least for the (!) dense stratum $S_d^0$.
\end{itemize}

The most favourable property regarding quotients is the following:

\begin{hope} For certain ``good" stratifications, one can relate (in a geometric way) a quotient $S_i//G_d$ to the structure of some $S_{d'}^0//G_{d'}$ for another dimension vector $d'$, or to some $S_{d'}^0(Q')//G_{d'}$ for another quiver $Q'$, or to some $(k^{n\times n})^m//{\rm GL}_n$.
\end{hope}

There are essentially two sources for such stratifications: on the one hand, ``geometric stratifications", obtained by specializing general concepts from Algebraic Geometry, Geometric Invariant Theory, or the theory of transformation groups, to the special situation of the action $G_d:R_d$. On the other hand, ``representation theoretic stratifications", which are constructed using representation theoretic concepts.\\[2ex]
In general, the ``geometric stratifications" have the advantage that, just from the very beginning, some good geometric properties are known. But the disadvantage usually is that the relation to the representation theory of $Q$ is difficult to get; in particular, their combinatorics is usually difficult to describe.\\[1ex]
Conversely, for the ``representation theoretic stratifications" one has a good control over the combinatorial aspects, or nice reduction techniques, etc., but it is then usually very difficult to establish any good geometric properties, especially concerning the existence (and description) of quotients.\\[1ex]
Some examples for the first type of stratifications include sheets \cite{BK}, the Hesselink stratification \cite{LB1,LB2,LB3}, or the Luna stratification of quotients \cite{LBP}.\\[1ex]
Representatively for the second type, which we will consider exclusively in the following, we have the stratification by decomposition type, the Harder-Narasimhan stratification, and the exhaustion of $R_d$ by subvarieties of fixed types of composition series. They will be discussed in the following section.

\subsection{Examples}\label{examples_strat}

As a first example, let us consider the stratification by decomposition type. Fix a dimension vector $d\in{\bf N}I$, and consider an arbitrary decomposition $d^*\, :\, d=d^1+\ldots+d^n$ into other dimension vectors. Define $S_{d^*}$ as the set of all representations isomorphic to a direct sum of indecomposables of dimension vectors $d^i$ for $i=1\ldots s$. From the point of view of the representation theory of finite dimensional algebras, this is the most natural stratification. It is also clear that we can easily characterize the non-empty strata by means of Kac's Theorem:
The stratum $S_{d^*}$ is non-empty if and only if each $d^i$ is a positive root.\\[1ex]
But this stratification usually has very bad geometric properties, as was already demonstrated in an example by H.~Kraft and C.~Riedtmann \cite{KRi}: Consider the quiver with set of vertices $\{i,j\}$, two arrows from $i$ to $j$, and one arrow from $j$ to $i$. Then, for the dimension vector $d=2i+j$, the dense stratum $S_{(d)}$ of indecomposables is not locally closed, but only constructible.\\[2ex]
As the next example, let us look at the Harder-Narasimhan filtration. 
The main result underlying its definition is the following:
\begin{proposition} Fix a slope function $\mu$ as in Section \ref{notation_moduli}. Then every representation $X$ possesses a unique filtration
$$0=X_0\subset X_1\subset\ldots\subset X_s=X$$
such that
\begin{itemize}
\item all subquotients $X_k/X_{k-1}$ are semistable,
\item $\mu(X_1/X_0)>\ldots>\mu(X_s/X_{s-1})$.
\end{itemize}
\end{proposition}

This allows us to make the following definition: suppose $d^*=(d^1,\ldots,d^s)$ is a decomposition of $d$, that is, $d=d^1+\ldots+d^s$. Then define the HN stratum $S_{d^*}$ as the set of all representations $X$ such that, in the Harder-Narasimhan filtration of $X$, we have $\dim X_k/X_{k-1}=d^k$ for all $k=1\ldots s$. Just from the definition, we see immediately that the open stratum $S_{(d)}=R^{ss}_d$ consists precisely of the semistable representations, so the results reviewed in Section \ref{notation_moduli} apply; in particular, the dense stratum admits an algebraic quotient.\\[1ex]
Moreover, one can relate the structure of an arbitrary HN stratum to the generic strata for other dimension vectors. More precisely, $S_{d^*}\simeq G_d\times^{P_{d^*}}V$ for a certain vector bundle $V\rightarrow\prod_{k=1}^sS_{(d^k)}$ and a certain parabolic subgroup $P_{d^*}\subset G_d$.\\[2ex]
Thus, we see that many of the general questions posed above allow an answer for the Harder-Narasimhan stratification. However, what is missing is an explicit description of the non-empty strata, or equivalently, an explicit criterion for the existence of semistables of a given dimension vector. The only known criteria are of a recursive nature (compare the discussion in \ref{general_rep}):

\begin{lemma} For $d\in{\bf N}I$, we have $R_d^{ss}\not=\emptyset$ if and only if there is no decomposition $d=d^1+\ldots+d^s$ such that $R_{d^k}^{ss}\not=\emptyset$ for all $k=1\ldots s$, $\mu(d^1)>\ldots>\mu(d^s)$ and $\langle d^k,d^l\rangle=0$ for all $k<l$.
\end{lemma}

As the third example, we formulate a proposal for a construction of a very general class of stratifications. The starting point is an arbitrary torsion theory $({\cal T},{\cal F})$ on $\bmod kQ$; that is, ${\rm Hom}(X,{\cal F})=0$ if and only if $X\in{\cal T}$, and ${\rm Hom}({\cal T},X)=0$ if and only if $X\in{\cal F}$, for all $X\in\bmod kQ$. Then one can define, for all pairs of dimension vectors $e\leq d$, a stratum $S_e$ in $R_d$ as the set of all $X\in R_d$ whose torsion part (with respect to $({\cal T},{\cal F})$) has dimension vector $e$. It seems that, using the same techniques as in the HN case, one can reduce arbitrary such strata to generic strata for smaller dimension vectors. Note that the HN stratification in fact fits into this general scheme: for any $\mu\in{\bf Q}$, define ${\cal T}_\mu$ as the set of all representations whose HN filtration have all subquotients of slope $\geq \mu$, and ${\cal F}_\mu$ as the set of all representations whose HN filtration have all subquotients of slope $<\mu$. It is then rather easy to see that $({\cal T}_\mu,{\cal F}_\mu)$ defines a torsion theory. This proposal suggests: 

\begin{hope} One can develop a general theory of stratifications of the $R_d$ induced by types of filtrations, torsion theories etc.~for $\bmod kQ$.
\end{hope}

\subsection{The composition monoid}\label{comp_monoid}

The last example to be discussed is of a quite different nature, since it is not a stratification, but an exhaustion of the representation varieties $R_d$. Furthermore, although its origin is representation theory, the tools to study its structure are inspired by quantum group theory. The following can be found in \cite{monoid}.\\[1ex]
Fix a dimension vector $d\in{\bf N}I$. For each word $\omega=(i_1\ldots i_n)$ in the alphabet $I$, being of weight $d$ (i.e.~such that $d=\sum_{k=1}^ni_k$), and any representation $X$, we define a composition series of type $\omega$ of $X$ to be a sequence
$$X=X_0\supset X_1\supset\ldots\supset X_n=0$$
such that $X_{k-1}/X_k\simeq E_{i_k}$ for $k=1\ldots n$, where $E_i$ denotes the simple representation corresponding to a vertex $i$ in $I$. Finally, we define ${\cal E}_{\omega}$ as the set of all representations $X\in R_d$ possessing a composition series of type $\omega$. The ${\cal E}_\omega$ are closed irreducible subvarieties of $R_d$.\\[1ex]
Denote the set of all words as above by $\Omega$. This set becomes a monoid via concatenation of words. We also define a partial ordering on $\Omega$ as follows: choose a total ordering on $I$ such that the existence of an arrow $i\rightarrow j$ implies $i<j$ in $I$. Then define the partial ordering on $\Omega$ to be generated by
$$\omega_1 ij\omega_2<\omega_1 ji\omega_2$$
whenever $i<j$ in $I$, and $\omega_1,\omega_2\in\Omega$. The following is easy to see:
\begin{lemma} If $\omega\leq\omega'$, then ${\cal E}_\omega\supset{\cal E}_{\omega'}$.
\end{lemma}
Thus, we see that the ${\cal E}_\omega$ for words $\omega$ of weight $d$ define an exhaustion of $R_d$: we always have the trivial stratum $\{0\}={\cal E}_{n^{d_n}\ldots 1^{d_1}}$, and the full representation variety itself $R_d={\cal E}_{1^{d_1}\ldots n^{d_n}}$, where we assume $I=\{1,\ldots,n\}$, compatible with the total ordering chosen above.\\[1ex]
The main question about this stratification is under which conditions on the words $\omega$, $\omega'$ do we have equality ${\cal E}_\omega={\cal E}_{\omega'}$ of the corresponding closed subvarieties of $R_d$. Certain degenerate variants of quantum groups enter naturally in this problem, which is described in the following.\\[2ex]
The basis for this construction is:
\begin{lemma}\label{lmot} Let ${\cal A}\subset R_d$ and ${\cal B}\subset
R_e$ 
be irreducible, closed subsets,
stable under the action of the group $G_d$
(resp.~$G_e$). Define ${\cal A}*{\cal B}$ as the set of all
representations 
$X\in R_{d+e}$ such that there exists an exact
sequence
$\ses{B}{X}{A}$ for some $A\in{\cal A}$, $B\in{\cal B}$. Then ${\cal
A}*{\cal 
B}$ is again irreducible, closed,
and stable under the group $G_{d+e}$.
\end{lemma}

This allows the following definition:

\begin{definition} The set ${\cal M}$ of irreducible, closed, 
$G_d$-stable subvarieties of the varieties $R_d$ for various $d\in{\bf
N}I$, 
together with the operation $*$ and the unit element $R_0$, is called
the 
extension monoid of $Q$.
\end{definition}

The monoid ${\cal M}$ seems to be too large to study its algebraic structure 
in general. Therefore, in analogy to the Hall 
algebra situation in Section \ref{notation_hall}, we introduce the composition monoid ${\cal C}$ as the 
submonoid of ${\cal M}$ generated by the elements $R_i$ for $i\in I$. It is then clear that each element of ${\cal C}$ is of the form ${\cal E}_\omega$ for 
some word $\omega$. The relation of ${\cal C}$ to quantum groups is the following:

\begin{theorem}\label{t13} The map $\eta$ sending $E_i$ to $R_i$ for $i\in I$ 
extends to a 
surjective monoid homomorphism $\eta:{\cal U}\rightarrow {\cal C}$, where the 
monoid ${\cal U}$ is defined by generators
$E_i$ for $i\in I$ and relations
$$E_i^{n+1}E_j=E_i^nE_jE_i,\;\;\; E_iE_j^{n+1}=E_jE_iE_j^n$$
if there is no arrow from $j$ to $i$ in $Q$, and $n$ equals the number of 
arrows from $i$ to $j$.
\end{theorem}

Note that the monoid ring ${\bf C}{\cal U}$ is related to the specialization at $v=0$ of a twisted form of ${\cal U}_v^+(\mathfrak{g}_Q)$ (see Section \ref{notation_hall}); see \cite{monoid} for details.\\[1ex]
It is easy to give examples of comparison maps $\eta$ which are not bijective. Therefore, it becomes neccessary to develop 
new techniques for the study of relations in ${\cal C}$. It turns out that 
A.~Schofield's analysis of generic properties of representations, reviewed in Section \ref{general_rep}, 
fits into this context, in that these results can be reformulated as algebraic 
properties of ${\cal C}$. Two results in this direction are the following :

\begin{theorem}~
\begin{enumerate}
\item For $d,e\in{\bf N}I$, the relation $R_d*R_e=R_{d+e}$ holds in ${\cal C}$ 
provided ${\rm ext}(e,d)=0$.
\item Each element of ${\cal C}$ can be written in the form 
$R_{d^1}*\ldots*R_{d^s}$, where each $d^k$ is a Schur root, and ${\rm ext}(d^{k+1},d^k)\not=0$ provided ${\rm 
ext}(d^k,d^{k+1})\not=0$, for all $k=1\ldots s-1$.
\end{enumerate}
\end{theorem}

\section{Moduli spaces}\label{moduli}

The most accessible objects currently available for the geometric study of isomorphism classes of quiver representations are the moduli spaces ${\cal M}_d^{ss}$ defined in Section \ref{notation_moduli}. Some representation theoretists feel uncomfortable with them, since their definition depends on the choice of stability $\Theta$. Proposals for more canonical constructions exist, but apparently not much has been proved about them.\\[1ex]
So the general principle in the following is to ignore this choice, and to ask for a study of ${\cal M}_d^{ss}$ for arbitrary $\Theta$. This has the advantage of giving us the well-developed machinery of D.~Mumford's GIT \cite{Mu} at hand.

\subsection{Cohomology of moduli spaces}\label{cohomology}

The starting point of \cite{hnsystem} is the wish to translate the
techniques 
of \cite{HN} from the context of moduli spaces of
vector
bundles on smooth projective curves to the moduli spaces of quivers. In \cite{HN}, G.~Harder and M.~S.~Narasimhan provide a
recursive
technique to compute Betti numbers of such moduli spaces, which involves 
counting certain classes of vector bundles over finite fields,
together with the Weil conjectures \cite{De}. In a general context of 
Geometric Invariant Theory quotients (but using the Hesselink stratification), this technique is developed in
\cite{Kir1}. One of the key results of \cite{hnsystem} is that the HN 
recursion makes
sense 
even in the non-commutative context of quantum
groups. Whereas the original philosophy of \cite{HN} seeks to express the number
of 
rational points of $R_d^{ss}$ over a finite field $k$ by a
recursive expression, the idea of \cite{hnsystem} is to find a similar 
expression directly for the characteristic function
$E_d^{ss}$ of $R_d^{ss}$:
$$E_d^{ss}(X)=\left\{\begin{array}{ccc}1&,&X\in R_d^{ss}\\ 
0&,&\mbox{otherwise},\end{array}\right.$$
viewed as an element of the Hall algebra $H_v(Q)$. The result is the following:
 
\begin{theorem}\label{t31} For each $d\in{\bf N}I$, the element
$E_d^{ss}$ is 
given by
$$E_d^{ss}=E_d-\sum v^{-\sum_{k<l}\langle 
d^k,d^l\rangle}E_{d^1}^{ss}*\ldots*E_{d^s}^{ss},$$\label{rec}
where the sum runs over all sequences $(d^1,\ldots,d^s)$ of dimension
vectors 
which sum up to $d$, and which fulfill the property
$\mu(d^1)<\ldots<\mu(d^s)$. In particular, if $\Theta$ is constant on
the 
support of $d$, that is, on the $i\in I$ such that
$d_i\not=0$, then $E_d^{ss}=E_d$.
 
The element $E_d^{ss}$ is given explicitely by the formula
$$E_d^{ss}=\sum(-1)^{s-1} v^{-\sum_{k<l}\langle 
d^k,d^l\rangle}E_{d^1}*\ldots*E_{d^s},$$\label{resrec}
where the sum runs over all sequences $(d^1,\ldots,d^s)$ of non-zero 
dimension vectors which sum up to $d$, and which fulfill
the property $\mu(d)<\mu(d^k)+\ldots+\mu(d^s)$ for all $k=2\ldots s$.
\end{theorem}
 
\begin{corollary} For all $d\in{\bf N}I$, we have
$E_d^{ss}\in C_v(Q)\simeq U^+_v$.
\end{corollary}
 
We can therefore consider the elements $E_d^{ss}$ as elements of $U^+_v$. 
Since the composition algebra $C_v(Q)\subset H_v(Q)$ is defined by
generators 
only, it is difficult to decide for a given function
of $H_v(Q)$ whether it belongs to $C_v(Q)$ or not. This has been the subject 
of several works. The above
corollary 
provides a large class of such elements.\\[2ex]
The algebraic use of these results lies in the construction of a large 
linearly independent subset, the Harder-Narasimhan system,
in $U^+_v$.
 
\begin{theorem} The set 
$$\{E_{d^1}^{ss}\cdot\ldots\cdot E_{d^s}^{ss}\, :\, 
\mu(d^1)<\ldots<\mu(d^s)\}\setminus\{0\}$$
is an orthogonal system in ${U}^+_v$.
\end{theorem}
 
Turning to the geometry of the moduli spaces, one can specialize the
formulas 
of Theorem \ref{t31} at the evaluation character
$f\mapsto|G_d|^{-1}\sum_Xf(X)$
of $H_v(Q)$ to get a formula (recursive, 
resp.~closed) for the fraction 
$|R_d^{ss}|/|G_d|$. Using methods of Geometric Invariant
Theory \cite{Mu} over finite fields \cite[Section 6]{hnsystem}, together with 
the Weil conjectures 
(\cite{De}, applied similarly to \cite{Kir1} and \cite{Got}), one arrives at 
the
following formula for Betti numbers of moduli spaces:
 
\begin{theorem}\label{bfa} Suppose that $\Theta(d)$ and $\dim d$ are
coprime. 
Then ${\cal M}_d^{ss}({\bf C})$ is a smooth projective
complex variety, and its Poincar\'e polynomial in cohomology with complex 
coefficients $\sum_{i\in{\bf Z}}\dim_{\bf C}H^i({\cal M}^{ss}_d({\bf C}),{\bf 
C})v^i$ is given by the formula
$$(v^2-1)^{1-\dim d}v^{-\sum_{i\in I}d_i(d_i-1)}\sum(-1)^{s-1}v^{2\sum_{k\leq 
l}\sum_{i\rightarrow j}d_i^kd_j^l}\prod_{k=1}^s\prod_{i\in
I}([d^k_i]!)^{-1},$$
where the sum runs over all tuples $(d^1,\ldots,d^s)$ of non-zero dimension 
vectors which sum up to $d$, and which fulfill
the property $\mu(d)<\mu(d^k)+\ldots+\mu(d^s)$ for all $k=2\ldots s$.
\end{theorem}
 
This allows us to give new cohomology formulas even in classical situations of 
Geometric Invariant 
Theory, like the quotients $({\bf P}^{m-1})^n_{ss}/{\rm PGL}(m)$ of 
\cite[Chapter 3]{Mu}. See \cite[Section 7]{hnsystem} for more details. The Harder-Narasimhan techniques of \cite{hnsystem} have also proved to be 
useful in \cite{CBV}.\\[2ex]
It would be very desirable to derive similar results on counting special types of representations, most importantly stable representations, or simple representations in the more general case of quivers with oriented cycles. One first has to define the right objects to count:
\begin{definition} A representation $X$ of ${\bf F}_qQ$ is called absolutely simple (resp.~absolutely stable) if it remains simple (resp.~stable) after base extension to an algebraic closure $\overline{{\bf F}_q}$.
\end{definition}
For both situations, one can derive an arithmetic theory of such representations analogous to the one for absolutely indecomposables in \cite{Kac1, KRi}.
\begin{hope}\label{counting} One can derive recursive formulas for the number of (resp.~the number of isoclasses of) absolutely simple (resp.~absolutely stable) representations.
\end{hope}

\subsection{Moduli for Kronecker modules}\label{kronecker}

As a particular example of the use of the techniques developed above, let us look at moduli spaces for Kronecker modules. As in Section \ref{loop}, we consider the quiver $K_m$ with set of vertices $\{i,j\}$, and $m$ arrows from $i$ to $j$. We choose the canonical stability $\Theta=i^*$, and consider indivisible dimension vectors only.\\[1ex]
The quivers $K_m$ have a very special property, in that the reflection functors induce identifications of the moduli spaces (this can be checked rather easily by making the semistability condition explicit).

\begin{question} Are there any compatibilities between reflection functors and semistability for other quivers?
\end{question}

Therefore, we can assume that the dimension vector $d=d_ii+d_jj$ lies in the fundamental domain \cite{Kac1}, thus satisfying $\frac{2}{m}d_i\leq d_j\leq\frac{m}{2}d_i$. Moreover, there is a natural duality which allows us to assume $d_i\leq d_j$. Here are some examples of Betti numbers $(\dim H^{2i}({\cal M}_d^{ss}({\bf C}),{\bf C}))_i$ for the moduli spaces ${\cal M}_d^{ss}$ of the quiver $K_3$ under the above restrictions. They are obtained via an algorithm which is proved in \cite{hnsystem} together with Theorem \ref{bfa} above:
 
\begin{description}

\item[$d=i+j$:] 1,1,1

\item[$d=2i+3j$:] 1,1,3,3,3,1,1

\item[$d=3i+4j$:] 1,1,3,5,8,10,12,10,8,5,3,1,1 

\item[$d=4i+5j$:] 1,1,3,5,10,14,23,30,41,46,51,46,41,30,23,14,10,5,3,1,1 

\item[$d=5i+6j$:] 1,1,3,5,10,16,27,39,60,83,114,146,184,214,239,246,\ldots

\item[$d=5i+7j$:] 1,1,3,5,10,16,28,43,68,98,142,190,251,306,361,393,410,\ldots 

\item[$d=6i+7j$:] 1,1,3,5,10,16,29,43,69,100,149,206,289,380,504,635,792, 

942,1102,1221,1316,1339,\ldots

\end{description}

What can already be seen from this description is that, asymptotically, the Poincar\'e polynomial looks like
$$1+q+3q^2+5q^3+10q^4+16q^5+29q^6+\ldots$$
A quick glance at the Online Encyclopedia of Integer Sequences \begin{verbatim} http://www.research.att.com/~njas/sequences/ \end{verbatim} tells us that the coefficient of $q^n$ is the number $p_n$ of two-row plane partitions of weight $n$, that is, tuples of integers
$$\lambda_1\geq\lambda_2\geq\ldots\geq 0,\; \; \mu_1\geq\mu_2\geq\ldots\geq 0$$
such that $\lambda_i\geq \mu_i$ for all $i$, and $\sum_i\lambda_i+\sum_i\mu_i=n$.

\begin{question} Is there a natural combinatorial parametrization of cohomology classes in $H^*({\cal M}_d^{ss})$ by (certain special kinds of) such two-row plane partitions?
\end{question}

The generating function $\sum_{n=0}^\infty p_nq^n$ of two-row partitions equals:
$$(1-q)\cdot\prod_{i=1}^\infty(1-q^i)^{-2},$$
which is essentially the square of the Dedekind $\eta$-function, which becomes (almost) a modular form if $q$ is replaced by $e^{2\pi iz}$ (see also Section \ref{noncomm_geom} for further such modularity properties).

\begin{question} What (if any) is the geometric or representation theoretic relevance of this numerology?
\end{question}

That it is not too absurd to ask such questions is motivated by similar asymptotic results of J.-M.~Drezet \cite{Dr}. The idea is to fix the dimension vector, and to increase the number of arrows. At first sight, this is not too reasonable from the point of view of representation theory. But note that in the example to be worked out below in Section \ref{exkron22}, one can see that the moduli spaces becomes more well-behaved if the number of arrows is increased. This is also familiar in Classical Invariant Theory \cite{Kr1}. Drezet's result roughly states that, as the number $m$ of arrows increases, the moduli space ${\cal M}_{di+ej}^{ss}({\bf C})$ apporaches the classifying space $B(({\rm GL}_d\times{\rm GL}_e)/{\bf C}^*)$, whose Poincar\'e polynomial in cohomology is given by $(1-q)\cdot\prod_{i=1}^d(1-q^i)^{-1}\cdot\prod_{i=1}^e(1-q^i)^{-1}$.

\subsection{An example of singular moduli spaces}\label{exkron22}

As a particular example of moduli spaces for divisible dimension vectors, we consider the situation of the generalized Kronecker quiver $K_m$ with stability $\Theta=i^*$ and the dimension vector $d=2i+2j$. The techniques we will use in the following to describe the moduli space ${\cal M}_d^{ss}$ essentially come from Classical Invariant Theory, namely the techniques of polarization and restitution (see \cite{Kr1}).\\[1ex]
Denote by $Q_m$ the space of quadratic forms on the $m$-dimensional vector space $k^m$. We have a map $f$ from $R_d$ to $Q_m$ by assigning to a point in $R_d$, that is, a tuple $A=(A^1,\ldots, A^m)$ of $2\times 2$-matrices, a quadratic form $f_A$ defined by
$$f_A(\lambda_1,\ldots,\lambda_m)=\det(\sum_{i=1}^m\lambda_iA^i).$$
The key fact about this map is the following:
\begin{lemma} The tuple of matrices $A$ belongs to $R_d^{ss}$ if and only if the associated quadratic form $f_A$ is non-zero.
\end{lemma}
Studying the behaviour of the group action of $G_d$ on $R_d$ with respect to the map $f$ is easy:
$$f_{gAh^{-1}}=\frac{\det g}{\det h}f_A.$$
Thus, we see that the map $f$ becomes $G_d$-invariant once we pass to the projectivization ${\bf P}Q_m$ of $Q_m$. Note that this passage is made possible precisely by the above Lemma! We thus have a map, again denoted by $f$,
$$f:{\cal M}_d^{ss}\rightarrow{\bf P}Q_m.$$
To study this map in more detail, we use the additional symmetry of the quiver $K_m$. Namely, the group ${\rm GL}_m$ acts on $R_d$ by mapping a tuple of matrices to an appropriate linear combination:
$$g\cdot(A^1,\ldots,A^m)=(\sum_{i=1}^mg_{1i}A^i,\ldots,\sum_{i=1}^mg_{mi}A^i).$$
We also have an obvious action of ${\rm GL}_m$ on $Q_m$, and thus also on ${\bf P}Q_m$, via base change. Again, we have a nice compatibility:
\begin{lemma} For all $A\in R_d$ and all $g\in{\rm GL}_m$, we have $f_{g\cdot A}=g\cdot f_A$.
\end{lemma}
This additional group action gives a tool for determining the image of $f$: since ${\rm GL}_m$ acts on ${\bf P}Q_m$ with only finitely many orbits ${\cal O}_r$, in bijection with the possible ranks $r=1,\ldots, m$ of quadratic forms, we just have to find out which ranks of the quadratic forms $f_A$ can occur.

\begin{lemma} The rank of $f_A$ is less than or equal to $\min(4,m)$.
\end{lemma}
\begin{proposition} All possible ranks appear.
\end{proposition}
This is easily proven by looking at the following tuple of matrices:
$$\left[\begin{array}{ll}1&0\\ 0&1\end{array}\right],\; 
\left[\begin{array}{ll}i&0\\ 0&-i\end{array}\right],\; 
\left[\begin{array}{ll}0&1\\ -1&0\end{array}\right],\; 
\left[\begin{array}{ll}0&\pm i\\ \pm i&0\end{array}\right].$$

One can also check that the map $f:{\cal M}_d^{ss}\rightarrow{\bf P}Q_m$ is bijective over the orbits ${\cal O}_r$ for $r=1,2,3$. But it is branched over ${\cal O}_4$, since there is an additional invariant
$$\det\left[\begin{array}{llll}A^1_{11}&A^2_{11}&A^3_{11}&A^4_{11}\\
A^1_{12}&A^2_{12}&A^3_{12}&A^4_{12}\\
A^1_{21}&A^2_{21}&A^3_{21}&A^4_{21}\\
A^1_{22}&A^2_{22}&A^3_{22}&A^4_{22}\end{array}\right].$$
But this is the only obstruction to injectivity of $f$, so we finally get:
\begin{theorem} The moduli space ${\cal M}_d^{ss}$ admits a finite morphism f of degree $2$ onto the projectivization ${\bf P}\overline{{\cal O}_4}$ of the space of quadratic forms on $k^m$ of rank $\leq 4$, branched precisely over the quadratic forms of rank $4$.
\end{theorem}
By direct inspection of the above $4$-tuple of matrices, one can also easily describe how the moduli space of stable representations fits into our picture:
\begin{proposition} The set of stable representations is the inverse image of the quadratic forms of rank $3$ or $4$.
\end{proposition}
We have the following result describing the moduli spaces for small $m$:
\begin{corollary} For $m=1$, the moduli space ${\cal M}_d^{ss}$ reduces to a point, and there are no stable representations. For $m=2$, it is isomorphic to ${\bf P}^5$, and there are still no stables. For $m=3$, it is isomorphic to ${\bf P}^9$, and the set of stables is isomorphic to the homogeneous space ${\rm GL}_3/{\rm O}_3$. 
\end{corollary}
This point is a good opportunity to discuss what an ``explicit description" of a moduli space should be.
Look at ${\cal M}={\cal M}_d^{ss}$ for $m=4$. The above Theorem tells us that ${\cal M}$ can be realized by
taking ${\bf P}Q_4\simeq {\bf P}^{10}$, and taking a two-fold covering over the open subset ${\bf P}{\cal O}_4$,
which can be defined by the non-vanishing of the determinant, that is, by one homogeneous equation.
This is quite a handy description of ${\cal M}$, but maybe not what one would a priori call explicit.
But imagine we want to ``coordinatize" ${\cal M}$; that is, find an explicit embedding in some projective space.
Then the most obvious construction embeds ${\cal M}$ into ${\bf P}^{55}$, subject to many defining equations
of the image.
Surely this ``explicit" description does not give any hint on the structure of ${\cal M}$ any more.\\[2ex]
As a final remark, it seems that the above moduli spaces can also be interpreted in the framework of the theory
of compactifications of symmetric varieties of \cite{DCP}.

\subsection{The cohomology rings of moduli spaces}\label{cohom_ring}

As a natural next step in the cohomological study of the moduli spaces ${\cal M}={\cal M}^{ss}_d(Q)$ for indivisible $d\in{\bf N}I$, one can ask for the ring structure of cohomology $H^*({\cal M})=\oplus_{n}H^n({\cal M})$ under the cup product. \\[1ex]
A general technique for the description of this ring structure for quotient varieties was developed by G.~Ellingsrud and S.~A.~Str\o mme in \cite{ES}. Their results do not apply immediately to the case of quiver moduli. Nevertheless, A.~King, in unpublished work, has shown that the setup of \cite{ES} has to be modified only slightly to apply to quiver moduli. Instead of describing the result in general, we concentrate on the case of moduli for the generalized Kronecker quiver $K_m$ to give the reader some impression of the types of results to be expected. This is contained in \cite[Section 6]{ES}.\\[2ex]
So let $Q$ be the quiver $K_m$ and consider the dimension vector $di+ej$ for coprime $d,e\in{\bf N}$. The stability notion is again $\Theta=i^*$. Choose integers $a,b\in{\bf Z}$ such that $ae+bd=1$.  Consider the polynomial ring $A={\bf C}[\beta_1,\ldots,\beta_e,\gamma_1,\ldots,\gamma_d]$. We have an action of the product $W=S_e\times S_d$ of symmetric groups on $A$ via $(\sigma,\tau)\beta_i=\beta_{\sigma(i)}$ and $(\sigma,\tau)\gamma_i=\gamma_{\tau(i)}$. The invariant subalgebra $A^W$ is freely generated by the elementary symmetric functions in the $\beta_i$ and in the $\gamma_i$, that is, by the $b_i$ for $i=1\ldots e$ and the $c_i$ for $i=1\ldots d$, where
$$b_i=\sum_{1\leq n_1<\ldots<n_i\leq e}\beta_{n_1}\cdot\ldots\cdot\beta_{n_i},$$
and similarly for the $c_i$.\\[1ex]
There is a map from $A$ to the invariant subalgebra $A^W$ defined by
$$p(r)=\frac{\sum_{(\sigma,\tau)\in W}{\rm sgn}\sigma\cdot{\rm sgn}\tau\cdot(\sigma,\tau)(r)}{\prod_{i<j}(\beta_i-\beta_j)\cdot\prod_{i<j}(\gamma_i-\gamma_j)}.$$
Define the ideal $I$ of $A$ to be generated by the elements
$$a\cdot\sum_{i=1}^e\beta_i+b\cdot\sum_{i=1}^d\gamma_i$$
and the
$$\prod_{j=1}^d\prod_{i=l}^e(\beta_{\sigma(i)}-\gamma_{\tau(j)})^m$$
for $(\sigma,\tau)\in W$ and $k=1\ldots m$, where $l$ is defined as the least integer $\geq\frac{k\cdot e}{d}$.
\begin{theorem} The cohomology ring $H^*({\cal M}_{di+ej}^{ss}(K_m),{\bf C})$ is isomorphic to the quotient of $A^W={\bf C}[b_1,\ldots,b_e,c_1,\ldots,c_d]$ by the image $p(I)$ of the ideal $I$ under the map $p:A\rightarrow A^W$.
\end{theorem}
As an exercise in understanding this description, the reader might work out the case $d=1$, in which case we know that the moduli space is isomorphic to the Grassmannian ${\rm Gr}_e^m$, whose cohomology is well-known (see \cite{Fu})\\[2ex]
By a result of A.~King and C.~Walter \cite{Ki2}, the cohomology ring is generated by the Chern classes of the so-called universal bundles. In the above description, the generators $b_i$, $c_i$ do indeed correspond to these Chern classes.\\[2ex]
In principle, the above description of the cohomology ring should allow us in particular to read off the Hilbert series of that ring; that is, the Poincar\'e polynomial of cohomology.

\begin{question} How can one relate the above result to the description of the Poincar\'e polynomial given by \ref{bfa}? 
\end{question}

In the case of Grassmannians (or of flag manifolds), the description of the cohomology ring links to the modern formulation of Schubert calculus \cite{Fu}. In this theory, one has two descriptions of cohomology, one by an explicit basis (the Schubert cycles), the other by generators and relations. 
The main goal is then to express the basis elements explicitly in terms of the generators.

\begin{hope} A similar theory can be developed for (some) smooth projective quiver moduli. More precisely, one can find explicit bases for the cohomology, and explicit formulas expressing these basis elements in terms of Chern classes of universal bundles.
\end{hope}

\subsection{Further analogies to vector bundles}\label{vb}

The above results, namely the description of the Betti numbers of moduli spaces and the description of cohomology rings, are well known in the theory of moduli of vector bundles. Some of the relevant literature is \cite{AB,HN,Kir1,Kir2,Kir3,Kir4,Kir5}. Therefore, it is natural to ask which further topics in this theory might have analogies in the theory of quiver moduli. This suggests many questions and directions for future work.

\begin{question} What is the ``right" cohomology theory to be considered for moduli spaces ${\cal M}_d^{ss}$ or the subvariety of stables ${\cal M}_d^s$ when $d$ is divisible?
\end{question}

Candidates are, for example, cohomology with compact support of ${\cal M}_d^s$, since this is a smooth, but not neccessarily projective, variety. Or one could study the global intersection cohomology of ${\cal M}_d^{ss}$, since this is a projective, but not neccessarily smooth, variety.

\begin{question} Can one construct smooth compactifications of the ${\cal M}_d^s$?
\end{question}

Note that ${\cal M}_d^{ss}$ is a compactification of ${\cal M}_d^s$, but not smooth.

\begin{question} Can one translate the methods of \cite{Kir2} to construct (partial) desingularizations of the ${\cal M}_d^{ss}$? If so, which representation theoretic objects (if any) do these varieties parametrize? For example, can they be interpreted as moduli for quiver representations with an additional structure of some type (like certain compatible filtrations, extension classes, etc.)?
\end{question}

The main working tool in the series of articles of F.~Kirwan is the Hesselink stratification. Since the recursive description of the Poincar\'e polynomial of cohomology of the moduli spaces of \cite{hnsystem} was possible with the aid of the Harder-Narasimhan stratification, which is much easier to control than the Hesselink stratification, one has the following question:

\begin{question} Is there a relation between the Harder-Narasimhan stratification and the Hesselink stratification?
\end{question}

Motivated by \cite{Kir4}, we have:

\begin{question} Can one describe defining relations between the Chern classes of universal bundles conceptually; that is, from certain properties of the bundles themselves?
\end{question}

\section{Noncommutative algebraic geometry}\label{noncomm_geom}

There are several versions of what could be noncommutative algebraic geometry. We consider in the following the proposal of M.~Kontsevich and A.~Rosenberg \cite{KR}. The idea is to view (certain) noncommutative rings $A$ as ``machines" producing infinitely many (commutative) varieties, which should be viewed as an approximation to the ``non-commutative geometry" of $A$.\\[1ex]
One candidate for such an infinite series of varieties produced from a noncommutative ring is the series of varieties of representations $\bmod_n(A)$ for $n\in{\bf N}$. The rings $A$ serving as good sources are the so-called formally smooth (or quasi-free) algebras. They give rise (under certain mild conditions) to a series of smooth varieties $\bmod_nA$. Note that all these algebras are hereditary, and thus, one should view them as noncommutative curves.\\[2ex]
One important class of formally smooth algebras are the path algebras of quivers $Q$. Instead of the varieties of $n$-dimensional modules $\bmod_nkQ$, we can also consider the series of representation varieties $R_d$ for $d\in{\bf N}I$ by \cite{Bo}. In this respect, this approach of noncommutative geometry is what we have considered all the time!\\[1ex]
The only new point of view is that one should always consider all the $R_d$ at the same time, and that one should consider them as representing only a small portion of the full noncommutative geometry of $kQ$ itself. In this respect, one might view the study of, for example, single moduli spaces ${\cal M}_d^{ss}(Q)$ as not really adapted to the present situation. What are the consequences of that point of view?
In the introduction to his book \cite{Nak}, H.~Nakajima asks whether certain ``formal series of varieties" $\sum_nX_nz^n$ for varieties $X_n$ make sense.\\[1ex]
In the present setup, one could, for example, ask for the meaning of the formal series $\sum_n{\cal M}_{nd}^{ss}(Q)$, where $d$ is an indivisible dimension vector. Or, more generally, for the formal series $\sum_{d:\mu(d)=\mu}{\cal M}_d^{ss}(Q)$ for some fixed slope $\mu\in{\bf Q}$. Possible questions in this direction are:

\begin{question} What is the nature of the generating series
$$\sum_{n=0}^\infty\chi({\cal M}_{nd}^{(s)s}(Q))z^n,$$
where $\chi$ denotes the Euler-characteristic for a suitable cohomology theory (see \ref{vb} for a discussion).
\end{question}

A conceptual approach to such a question is provided by:

\begin{hope} The direct sum of cohomologies
$$\bigoplus_{n=0}^\infty H^*({\cal M}_{nd}^{(s)s}(Q))$$
(again for a suitable cohomology theory) can be equipped with an additional structure, for example, the structure of a representation of some (Lie) algebra.
\end{hope}

The above questions should not be taken too literally; there are many possible variations one can think of, some of which are described briefly in the following examples.\\[2ex]
The first example concerns work of A.~Klyachko \cite{Kly} on moduli of vector bundles on the projective plane. He essentially reduces the study of the cohomology of such moduli spaces to moduli spaces of representations of a star quiver with three arms, {\it for various stabilities $\Theta$}. In the case of vector bundles of rank $2$, one ends up with the study of representations of a quiver of type $D_4$, the dimension vector being the highest root, which is of course a well-known situation. The result is that the Euler characteristics of these moduli spaces of vector bundles are given by the Hurwitz function counting equivalence classes of integral binary quadratic forms. These results have been reinterpreted by C.Vafa and E.~Witten \cite{VW} to the extent that the generating function of Euler characteristics of these moduli spaces, for various second Chern classes of the bundles, is related to Eisenstein series, and thus to modular forms.\\[2ex]
The second example concerns Hilbert schemes of points in the affine plane, see e.~g.~\cite{Nak,Got}. By results of \cite{Got}, the generating functions of Euler characteristics are related to the Dedekind $\eta$-function. These results were interpreted conceptually by H.~Nakajima \cite{Nak}, by showing that the direct sum of cohomology spaces of Hilbert schemes of various numbers of points carries the structure of an irreducible representation of (a variant of) an infinite Heisenberg algebra. These results are related to quivers since it can be shown that these Hilbert schemes are essentially moduli spaces for representations of the deformed preprojective algebras associated to the quiver with set of vertices $\{i,j\}$, an arrow from $i$ to $j$, and a loop at $j$.\\[2ex]
Finally, the paper \cite{Ko} of M.~Kontsevich should be mentioned. It contains many fascinating speculations and proposals for noncommutative algebraic geometry (in the framework of \cite{KR}), some of which have direct connections to quiver representations. One of his proposals is to form generating functions of numbers of points of the module varieties $\bmod_n(A)$ of formally smooth algebras (so, for example $A=kQ$), and to study their modular properties.

\end{document}